\documentclass{amsart}
\usepackage{amsmath}
\usepackage{amssymb} 
\usepackage{enumerate}

\usepackage{graphicx}

\newtheorem{theorem}{Theorem}
\newtheorem{proposition}[theorem]{Proposition}
\newtheorem{lemma}[theorem]{Lemma}
\newtheorem{claim}[theorem]{Claim}

\begin{document}

\title{
When is a band-connected sum equal to the connected sum?${}^*$
}

\footnotetext[1]{
This question has been posted on
Mathematics Stack Exchange:
https://math.stackexchange.com/questions/2503991/when-is-a-band-connected-sum-equal-to-a-connected-sum-of-knots.
}

\author[K. Miyazaki]{Katura Miyazaki}
\address{Faculty of Engineering, Tokyo Denki University, 5 Senju Asahi-cho, Adachi-ku, Tokyo 120--8551, Japan}
\email{miyazaki@cck.dendai.ac.jp}
\date{}


\begin{abstract}
We show that
a band-connected sum of knots $K_0$ and $K_1$ along a band $b$
is equal to the connected sum $K_0\# K_1$ if and only if
  $b$ is a trivial band.
\end{abstract}

\maketitle

\section{Introduction}
\label{section:introduction}
Let $K_0 \cup K_1$ be a $2$-component split link in $S^3$,
and $b: I\times I\to S^3$ an embedding satisfying
$b(I\times I) \cap K_i = b(\{i\} \times I )$ for $i=0,1$.
The image $b(I\times I)$ is called a band connecting $K_0$ and $K_1$;
by abuse of notation, $b$ also denotes $b(I\times I)$.
Then the knot $K_0\cup K_1 \cup \partial b -\mathrm{int}( (K_0 \cup K_1) \cap b)$ is called the \textit{band-connected sum} of $K_0$ and $K_1$ along $b$, and 
denoted by $K_0 \#_b K_1$,
or $K_b$ for short.

A core of a band $b$ is $b(I\times \{*\})$, and a cocore of $b$
is $b(\{*\}\times I)$, where $*\in \mathrm{int}{I}$.
A band $b$ is \textit{trivial} if there is a splitting sphere for $K_0\cup K_1$ meeting $b$ in a cocore of $b$.
Obviously, if $b$ is a trivial band, then $K_0 \#_b K_1 \cong K_0\#K_1$.
We prove  the converse.

\begin{theorem}
\label{th}
If a band-connected sum $K_0 \#_b K_1$ is equal to the connected sum $K_0 \# K_1$, then $b$ is a trivial band.
\end{theorem}

Regarding the problem when a band-connected sum is composite, 
Eudave-Mu\~noz proved the following.

\begin{theorem}[\cite{EM}]
\label{EMth}
If $K_b$ is a composite knot,
then there is a decomposing sphere for $K_b$ which is either
disjoint from $b$ or intersects $b$ in a core of $b$.
\end{theorem}

To prove Theorem~\ref{th}
we use Theorem~\ref{EMth} and the fact that
if $K_b \cong K_0 \# K_1$, $K_b$ bounds a Seifert surface
of type~1a (Lemma~\ref{lemma}), which is  defined below.
Any band-connected sum $K_b$ bounds a compact connected surface $S$ containing $b$.
Such a surface $S$ is called \textit{type~1}
if $S-b$ is disconnected,
and \textit{type~2} otherwise.
In other words,
a type~1 surface for $K_0\#_b K_1$ is a compact connected,
possibly nonorientable, surface that is the union of $b$ and
two connected surfaces bounded by $K_i$ $(i=0,1)$.
A surface $S$ of type~1 is called \textit{type~1a} if there is 
a splitting sphere for $K_0\cup K_1$ that is disjoint from $S -b$,
and\textit{ type~1b} otherwise.

For a composite band-connected sum bounding a surface of type~1a,
we study the configuration of a splitting sphere, a decomposing sphere,
and a type~1a surface,
and obtain Theorem~\ref{mainth} below.
Theorem~\ref{th} follows from Theorem~\ref{mainth}.

\begin{theorem}
\label{mainth}
If $K_b$ is a composite knot and bounds a surface of type~$1$a,
then there is a decomposing sphere disjoint from $b$.
\end{theorem}

As another corollary to Theorem~\ref{mainth},
we give sufficient conditions for band-connected sums to be prime.

\begin{proposition}
\label{prop}
Assume that $K_b = K_0 \#_b K_1$ bounds
a surface $S_0 \cup b \cup S_1$ of type~$1$a, where
$S_i$ is a surface bounded by $K_i$ $(i=0,1)$.
Then $K_b$ is prime if $b$ is a nontrivial band
 and either $(1)$ or $(2)$ below holds.
\begin{enumerate}[$(1)$]
\item $K_i$ is a trivial knot for $i=0,1$.
\item $K_i$ is prime and $S_i$ is an incompressible Seifert surface for 
$i=0,1$.
\end{enumerate}
\end{proposition}

\section{Proofs}
\label{section:proofs}
\begin{lemma}
\label{lemma}
If $K_b \cong K_0\# K_1$, then $K_b$ bounds a Seifert surface of type~$1$a.
\end{lemma}

\noindent
\textit{Proof of Lemma}~\ref{lemma}.
The assumption implies $g(K_b) = g(K_0) +g(K_1)$,
where $g(\cdot)$ denotes the genus of a knot.
Then, by \cite{Ga, Sch} 
$K_b$ bounds a Seifert surface $S$ of type~1 that is the union of
$b$ and minimal genus Seifert surfaces $S_i$ for $K_i$ $(i=0,1)$.
Let $P$ be a splitting sphere for $K_0 \cup K_1$ such that
$|P \cap (S -b)|$ is minimal among all splitting spheres.
Suppose  $P \cap (S-b) \ne \emptyset$ for a contradiction;
then $P\cap (S-b)$ consists of simple closed curves.
Let $\Delta_1$ be an innermost one
among disks in $P$ bounded by components of $P\cap(S-b)$.
Without loss of generality $\partial \Delta_1 \subset S_0$.
If $\partial \Delta_1$ is essential in $S_0$, then
surgery of $S_0$ along $\Delta_1$ yields a Seifert surface for $K_0$
with fewer genus than $S_0$ and possibly a closed surface.
Since $S_0$ has minimal genus, we see $\partial \Delta_1$ bounds a disk
$\Delta_2$ in $S_0$.

Each component of $\Delta_2 \cap P$ bounds a disk in $\Delta_2$,
and let $\Delta_3 (\subset S_0)$ be an innermost one among all
such disks in $\Delta_2$;
note $\Delta_3 \cap P = \partial \Delta_3$.
Let $B$ be the $3$-ball satisfying $\partial B = P$ and $\Delta_3\subset B$,
and $B'$ the closure of a component of $B -\Delta_3$ containing $K_0$
or $K_1$.
Shrinking $B'$ slightly, we obtain a splitting sphere $\partial B'$ with 
$|\partial B' \cap (S-b)| < |P \cap (S-b)|$.
This contradicts the minimality of $|P\cap (S-b)|$.
It follows $P\cap (S-b) =\emptyset$, so that $S$ is type~1a
\hfill $\square$

\medskip
\noindent
\textit{Proof of Theorem}~\ref{mainth}.
By Theorem~\ref{EMth} we may assume that there is a decomposing
sphere $Q$ for $K_b= K_0\#_b K_1$ intersecting $b$ in a core of $b$.
Take a splitting sphere $P$ for $K_0 \cup K_1$,
a decomposing sphere $Q$ for $K_b$, and a type~1a surface $S$
bounded by $K_b$ that satisfy conditions~(1)---(4) below.
\begin{enumerate}[(1)]

\item $P, Q, S$ are in general position, i.e. any two of these
are in general position, and the intersection of any two of these
and the rest are in general position.

\item $P \cap (S -b) = \emptyset$.

\item $Q \cap b$ is a core of $b$, and $P\cap b$ consists of
cocores of $b$.

\item $|Q \cap S| + |P \cap Q| + |P \cap Q \cap S|$ is minimal
among all $P, Q, S$ satisfying (1)---(3) above.
\end{enumerate}

We first study the configuration of $P \cap Q$ and $Q \cap S$ on $Q$.
Note that $Q \cap S$ consists of an arc connecting the two points $K_b \cap Q$ and possibly simple closed curves.
Let $\alpha$ be the arc component of $Q \cap S$, and
set $c = Q\cap b$, a core of $b$.
Since $c$ connects  $K_0$ and $K_1$, and $S -b$ is disconnected,
$c$ is contained in $\alpha$.
Note also that $P\cap Q$, consisting of simple closed curves,
is disjoint from any circle component of $Q\cap S$
because $P\cap (S-b) =\emptyset$.
For the same reason $P\cap Q$ intersects $\alpha$ only in $c$.
See Figure~1.

\begin{figure}[h!]
\includegraphics[clip]{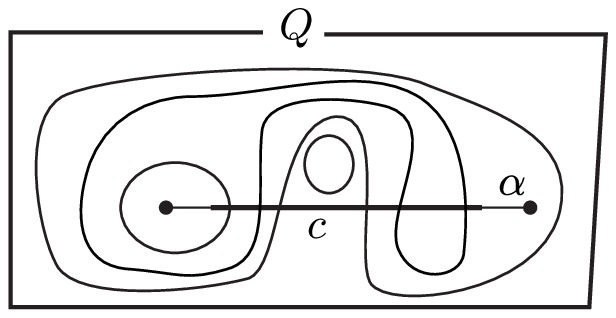}
\caption{$P\cap Q$ consists of simple closed curves.
\newline
$Q\cap S$ consists of the arc $\alpha$ and simple closed curves.}
\end{figure}

Assume for a contradiction that
either $Q \cap S$ or $P\cap Q$ has
a circle component disjoint from $\alpha$.
Each such component bounds a disk in $Q$ disjoint from $\alpha$.
Let $\Delta$ be an innermost one among all these disks.
If $\partial \Delta$ is a component of $Q \cap S$,
surger $S$ along $\Delta$ to yield a bounded surface $S'$ and
possibly a closed surface.
Since $\Delta \cap P = \emptyset$,
we see that $P \cap (S' -b) = \emptyset$ and
$|Q \cap S'| < |Q \cap S|$.
This contradicts the minimality assumption~(4).

If $\partial \Delta$ is a component of $P \cap Q$,
surger $P$ along $\Delta$ to yield two spheres,
one of which is  a splitting sphere $P'$.
Since $\Delta \cap S =\emptyset$,
we see that $P' \cap (S -b) = \emptyset$ and
$|P' \cap Q| < |P\cap Q|$.
This contradicts the minimality assumption.
We thus obtain the following.

\begin{claim}
\label{no disjoint circles}
$Q\cap S = \alpha$, and each component of $P\cap Q$ intersects
$\alpha$. \hfill $\square$
\end{claim}

Assume for a contradiction that
there is a component of $P\cap Q$ intersecting $\alpha$
more than once.
Each such component bounds a disk in $Q$ containing 
at most one endpoint of $\alpha$.
Let $\Delta_1$ be an innermost one among all such disks.
Then, $\Delta_1 \cap \alpha$ consists
of at most one non-properly embedded arc and at least one
properly embedded arc (Figure~2(1)).
Each properly embedded arc of $\Delta_1 \cap \alpha$
and a subarc of $\partial \Delta_1$ cobound a disk in $\Delta_1$
containing no endpoint of $\alpha$.
Let $\Delta_2$ be an innermost one among all such disks;
note $\Delta_2 \cap (S-b) =\emptyset$.
We can isotop $b$ along $\Delta_2$ to reduce $|P\cap Q\cap S|$
(Figure~2(2)), a contradiction to the minimality.
Thus Claim~\ref{intersectonce} below is proved.

\begin{figure}[h!]
\includegraphics[clip]{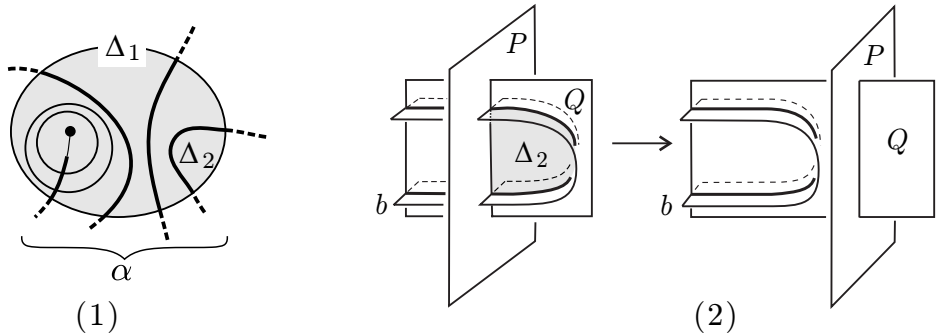}
\caption{Circles in (1) are components of $P\cap Q$.}
\end{figure}

\begin{claim}
\label{intersectonce}
Each component of $P\cap Q$ intersects $\alpha$ in one point of the core $c$.
\hfill $\square$
\end{claim}

\begin{figure}[h!]
\includegraphics[clip]{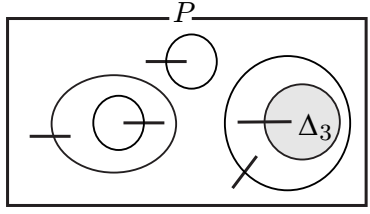}
\caption{Circles are the components of $P\cap Q$.
Arcs are the components of $P\cap b$, i.e. cocores of the band.}
\end{figure}

Let us take a look at $P\cap Q$ on the splitting sphere $P$.
Let $\Delta_3$ be an innermost one among all disks in $P$
bounded by components of $P\cap Q$; see Figure~3.
Surger the decomposing sphere $Q$ along $\Delta_3$ to obtain
two spheres, at least one of which is a decomposing sphere for $K_b$.
Let $Q'$ be a decomposing sphere obtained from $Q$ after the surgery.
Without loss of generality $Q'$ is disjoint from $K_0$
and intersects $K_1$ in two points.
Let $\Delta_4$ be the closure of the component of $b - (b\cap Q')$
disjoint from $K_0$ (Figure~4).
Sliding $Q'$ along the disk $\Delta_4$,
we obtain a decomposing sphere for $K_b$ disjoint from $b$ as claimed.
This completes the proof of Theorem~\ref{mainth}.
\hfill $\square$

\begin{figure}[h!]
\includegraphics[clip]{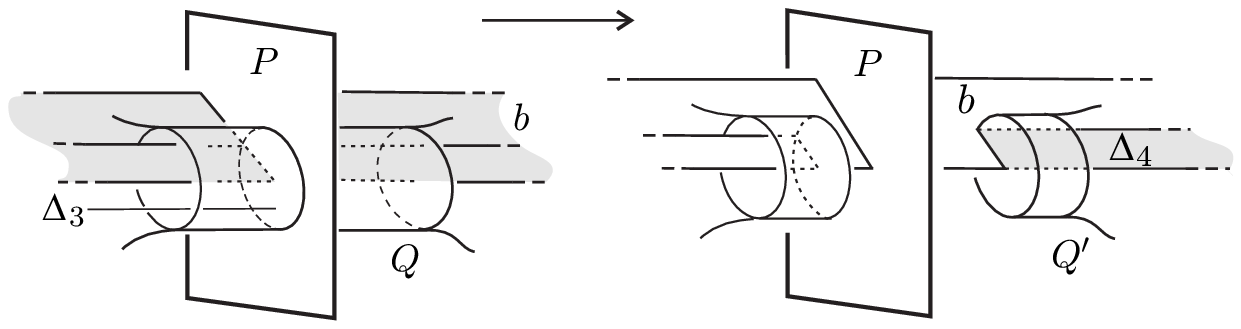}
\caption{}
\end{figure}

\medskip\noindent
\textit{Proof of Theorem} \ref{th}.
By Lemma \ref{lemma}
$K_b$ bounds a Seifert surface that is the union of minimal genus
Seifert surfaces $S_i$ for $K_i$ $(i=0, 1)$ and $b$.
We proceed by induction on the product of genera $g(K_0)g(K_1)$.

\textit{Case $1$. $g(K_0)g(K_1) =0$.}
\\
Since $K_0$ or $K_1$ is a trivial knot,
$S_0$ or $S_1$ is a disk.
The boundary of a regular neighborhood of the disk
 is a splitting sphere for $K_0\cup K_1$
intersecting $b$ in a single cocore of $b$.
It follows that  $b$ is a trivial band.

\textit{Case $2$. $g(K_0)g(K_1)\ge 1$.}
\\
By Theorem \ref{mainth} there is a decomposing sphere $Q$ for $K_b$
disjoint from $b$.
Without loss of generality $Q$ intersects $K_1$ in two points.
Let $B$ be the $3$-ball which is bounded by $Q$ and
disjoint from $K_0$.
Then,  $B\cap K_1$ is an arc knotted in $B$.
Let $K'_1$ be the knot obtained from $K_1$ by replacing
$(B, B\cap K_1)$ with $(D^2\times I, \{0\}\times I)$.
We can regard $b$ as a band connecting $K_0$ and $K'_1$.
If $b$ is a trivial band for $K_0\#_b K'_1$,
then it is trivial for $K_0\#_b K_1$.
Since $g(K_0)g(K'_1) < g(K_0)g(K_1)$,
Theorem~\ref{th} is proved inductively.
\hfill $\square$

\medskip\noindent
\textit{Proof of Proposition} \ref{prop}.
Assume for a contradiction that $K_b$ is not prime,
but satisfies either condition~(1) or (2) in Proposition~\ref{prop}.
If $K_b$ is a trivial knot, then $K_0$ and $K_1$ are trivial knots
and especially $b$ is a trivial band (\cite{Sch2}),
a contradiction to the assumption of Proposition~\ref{prop}.
It follows that $K_b$ is a composite knot.

Then, following the arguments in Case~2 in the proof of
Theorem~\ref{th},
without loss of generality we can take a $3$-ball $B$ such that
$B\cap (K_0 \cup b) =\emptyset$ but
 $B\cap K_1$ is an arc knotted in $B$.
This implies that $K_1$ is not a trivial knot,
which violates condition~(1).
It follows that $K_b$ satisfies condition~(2).

Let $K'_1$ be $K_1$ with $(B, B\cap K_1)$ replaced by
$(D^2\times I, \{0\}\times I)$.
Since $K_1$ is prime and $B\cap K_1$ is a knotted arc in $B$,
we see $K'_1$ is a trivial knot.
The Seifert surface $S_1$ for $K_1$ intersects $\partial B$ in
an arc connecting the two points $\partial B \cap K_1$ and
possibly simple closed curves.
However, using the incompressibility of $S_1$, we can isotop $S_1$
without meeting $S_0\cup b$ so that $S_1\cap \partial B$ is just one arc.
The arc $S_1\cap \partial B$ in $\partial B=
\partial (D^2\times I)$ and $\{0\}\times I$ cobound a disk
in $D^2\times I$.
Then, the union $S'$ of the disk and $S_1-\mathrm{int}B$ is
an incompressible Seifert surface for the trivial knot $K'_1$.
It follows that $S'$ and thus $S_1-\mathrm{int}B$ are disks.
Isotopy of $\partial B$ along the disk $S_1 -\mathrm{int}B$
yields a splitting sphere for $K_0\cup K_1$ intersecting $b$
in a cocore of $b$.
This contradicts the assumption that $b$ is a nontrivial band.
\hfill $\square$

\end{document}